\documentclass[10pt,a4paper,english]{article}
\usepackage[T1]{fontenc}
\usepackage[latin9]{inputenc}
\usepackage{textcomp}
\usepackage{graphicx}

\makeatletter

\pdfpageheight\paperheight
\pdfpagewidth\paperwidth

\usepackage{babel}

\makeatother

\usepackage{babel}
\begin{document}

\title{Rotationally Symmetric Tilings with Convex Pentagons and Hexagons }

\author{Bernhard Klaassen }

\date{Fraunhofer Institute SCAI\\
 Schloss Birlinghoven \\
 53754 Sankt Augustin, Germany\\
 bernhard.klaassen@scai.fraunhofer.de\\
$ $\\
Preprint to appear in Elemente der Mathmatik, Vol. 71 (2016)}
\maketitle
\begin{abstract}
In contrast to many known results concerning periodic tilings of the
Euclidean plane with pentagons, here tilings with rotational symmetry
are investigated. A certain class of convex pentagons is introduced.
It can be shown that for any given symmetry type $\mathbf{C}_{n}$
or $\mathbf{D}_{n}$ there exists a monohedral tiling generated by
a pentagon from this class. For $n>1$ each of these tilings is also
a spiral tiling with $n$ arms. As a byproduct it follows that the
same holds for convex hexagons. 
\end{abstract}

\section{Introduction}

Recent news in 2015 \cite{key-6} made a nearly hundred-years-old
problem popular again: the complete characterization of all convex
pentagons that tile the Euclidean plane monohedrally (i.e., filling
the plane completely without gaps or overlaps by using only congruent
copies of one special tile). There is no hint that this characterization
can be fully done in the near future, so only a systematic catalog
of the existing tilings can be given (e.g., see \cite{key-1} or \cite{key-2}).
But why are pentagons so interesting mathematically? The reason is
simple: For all other convex polygons the question is settled. So,
the pentagons remained as an open problem.

In most papers about this subject the emphasis is put on periodic
tilings. One reason might be the fact that all known families of pentagons
tiling the plane are able to generate periodic tilings but only few
of them also deliver rotationally symmetric ones (e.g., \cite{key-3}
or at the end of \cite{key-4}). It is good practice to denote the
type of symmetry in the catalogs of tilings. But we also could view
the problem from the other direction: For which types of symmetry
is a monohedral tiling with convex pentagons possible? We will answer
this question positively for all types of rotational symmetry. In
our case the Schoenflies notation for symmetry is the easiest one,
which uses $\mathbf{C}_{n}$ for $n$-fold rotational symmetry without
reflection and $\mathbf{D}_{n}$ for additional $n$ axes of reflection
(e.g., see fig. 1.1). 

\includegraphics[width=10cm]{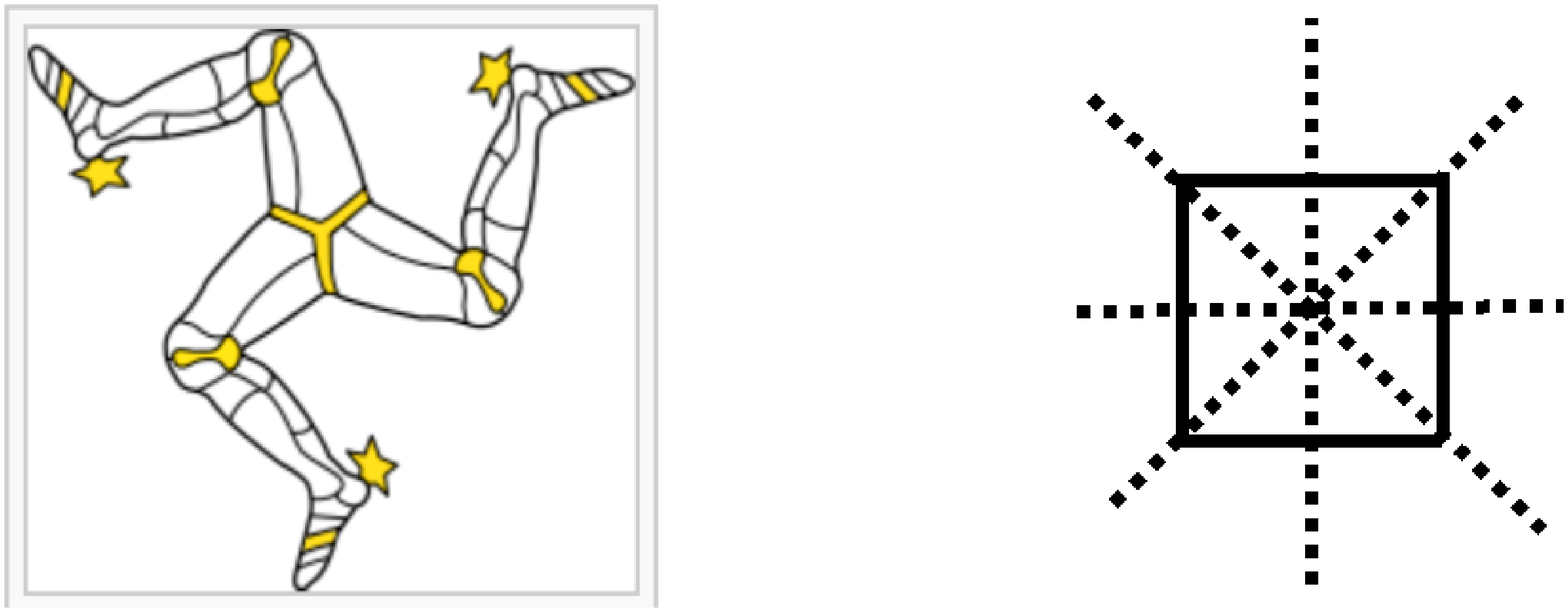} 
\begin{description}
\item [{Figure}] 1.1 \,\,\,Examples for $\mathbf{C}_{3}$ (taken from
\cite{key-5}) and $\mathbf{D}_{4}$ symmetry 
\end{description}
Regarding the periodic case, only four types of rotational symmetry
are possible: $n=2$, 3, 4 or 6. It is interesting to observe that
all of these cases occur in the catalog of known periodic pentagon
tilings. For the non-periodic case, which is studied here, we will
see that we can go further. From related work (e.\,g. \cite{key-11},
\cite{key-12}) tilings are known with arbitrary rotational symmetry
but they are generated by rhombs or triangles as prototiles. Here
we define a class of pentagons, from which one can generate the proposed
symmetry types. The following properties define the class:

\subparagraph*{Property 1:}

Let $A,\,B,\,C,\,D,\,E$ be the inner angles in anti-clockwise order
of a pentagon $P$ and $a,\,b,\,c,\,d,\,e$ the edges with $a$ ending
at corner $A$ and so forth. $P$ is regarded to have property 1 if
all inner angles $<180{^{\circ}}$ and $|b|=|c|=|a|+|d|$ and $D+E=180{^{\circ}}$
and ( $C|360{^{\circ}}\bigvee B|360{^{\circ}}$) where $|a|$ denotes
the length of side $a$. (See fig. 2.1 for an example.)

\ 

Obviously, this pentagon class is a subclass of ``type 1\textquotedblright{}
from \cite{key-2} and has three degrees of freedom: The angles $B$,
$C$ and $D$ can be chosen within certain ranges and constraints
but independently from each other. The other angles and proportions
are dependent from these three.

\section{Results}

Without loss of generality we assume for simplicity that angle $B$
divides 360\textdegree . (If $B$ and $C$ are interchanged, the following
theorem also holds with equivalent proof.)

\subparagraph*{Theorem:}

For a given natural number $n>2$ any property-1-pentagon with $B=360{^{\circ}}/n$
tiles the Euclidean plane with $n$-fold rotational symmetry. If $C=180{^{\circ}}-B/2$
and $D=90{^{\circ}}$, the plane can be tiled with type $\mathbf{D}_{n}$,
otherwise with $\mathbf{C}_{n}$. For $n=2$ or $n=1$ such a tiling
is possible with any $B<180{^{\circ}}$, especially with property-1-pentagons.

\subparagraph{Proof:}

We start with the case $n>2$. First, we should note that each pentagon
$P$ with property 1 has parallel sides $d$ and $a$. So we can take
a copy of $P$ and glue both together at side $e$ to form a hexagon,
in most cases not a regular one, but equilateral and with parallel
sides.

\,

\includegraphics[width=5cm]{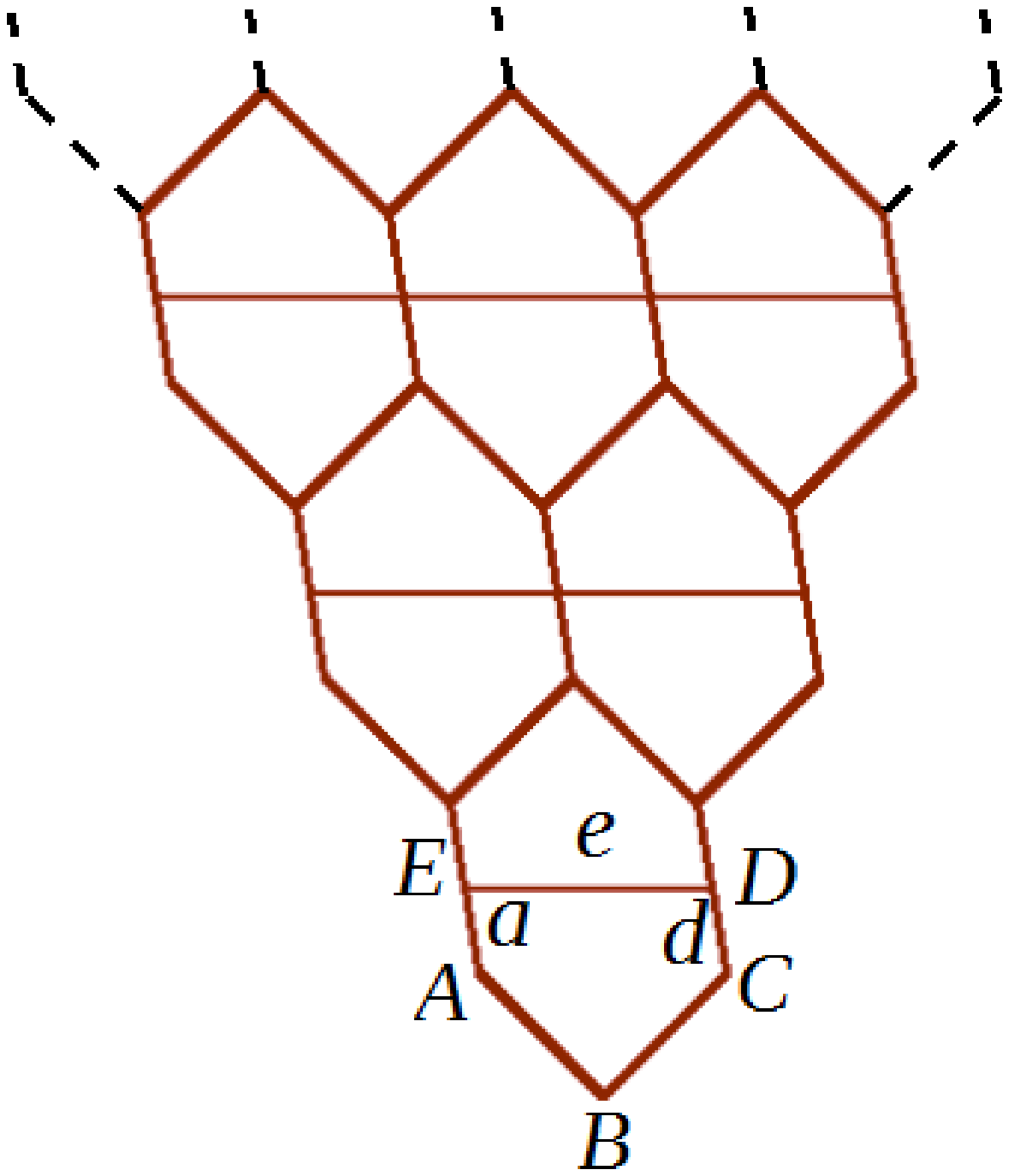} 
\begin{description}
\item [{Figure}] 2.1 \,\,\,A sector formed by congruent hexagons 
\end{description}
Putting copies of these hexagons together in form of growing rows,
one can fill a kind of sector of the plane. The outer shape of this
sector is fully characterized by the angles, since the length of the
line segments is always $|b|$. Viewing this sector in upright position
as in figure 2.1, we see that the left border is a zig-zag-line with
alternating inner angles $A$ and $180{^{\circ}}-A$, the right border
has angles $C$ and $180{^{\circ}}-C$ alternating.

Now we take a copy of this sector and connect it to the first one
as shown in figure 2.2. Between both sectors a gap remains which has
the following characterization: The innermost point has outer angle
$A+C$, which is $540{^{\circ}}-B-D-E=360{^{\circ}}-B$, so the inner
angle of the gap is $B$. The left border is an equilateral zig-zag
with alternating angles $C$ and $180{^{\circ}}-C$, the right border
has angles $A$ and $180{^{\circ}}-A$ alternating.

\,

\includegraphics[width=7cm]{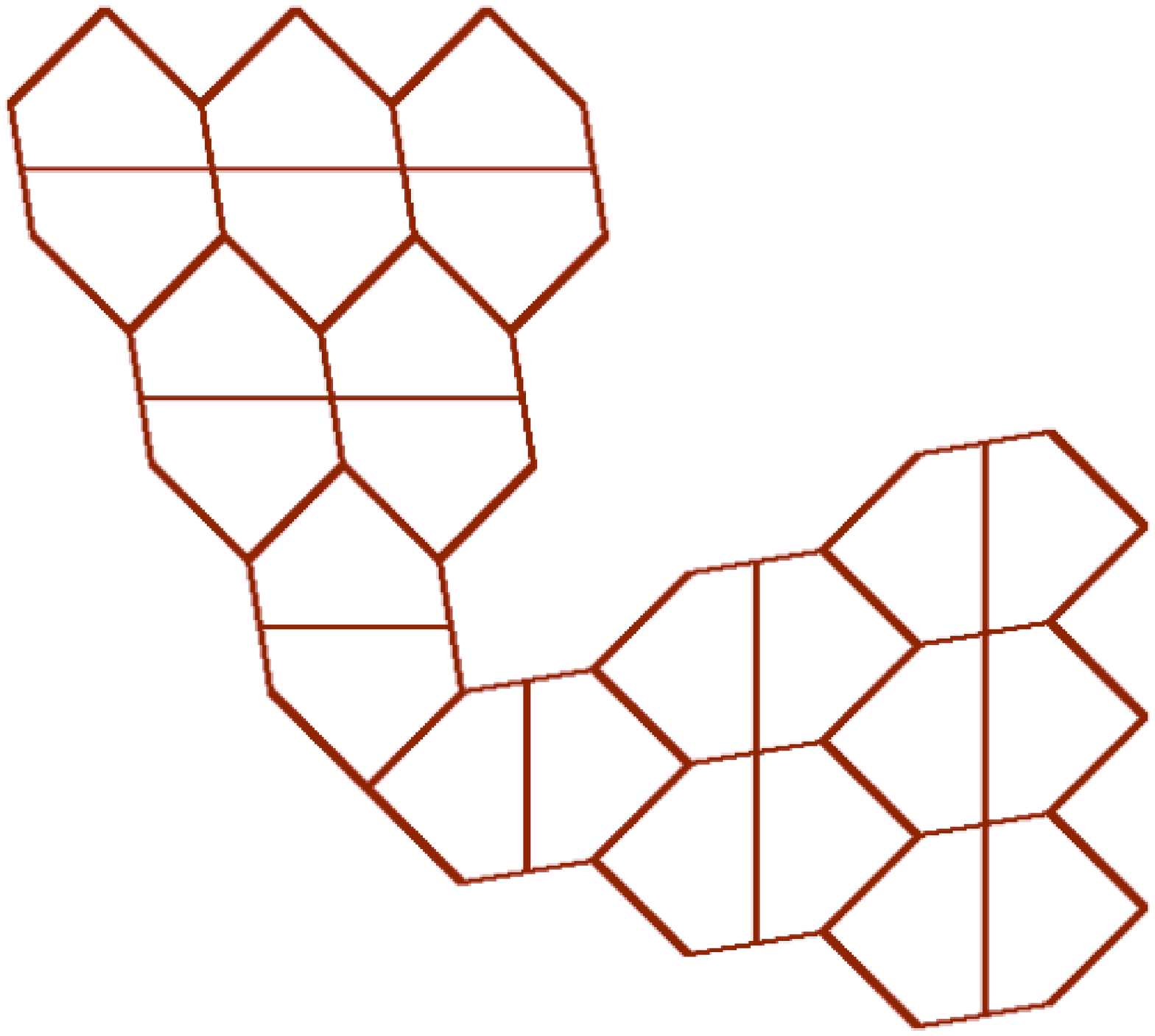} 
\begin{description}
\item [{Figure}] 2.2 \,\,\,Two sectors with a gap 
\end{description}
It follows that the original sector from figure 2.1 perfectly fits
into the gap after a reflection operation at the vertical axis and
a suitable rotation. Such a reflection plus rotation is allowed, since
all parts within a sector remain congruent under this operation.

\,

\includegraphics[width=7cm]{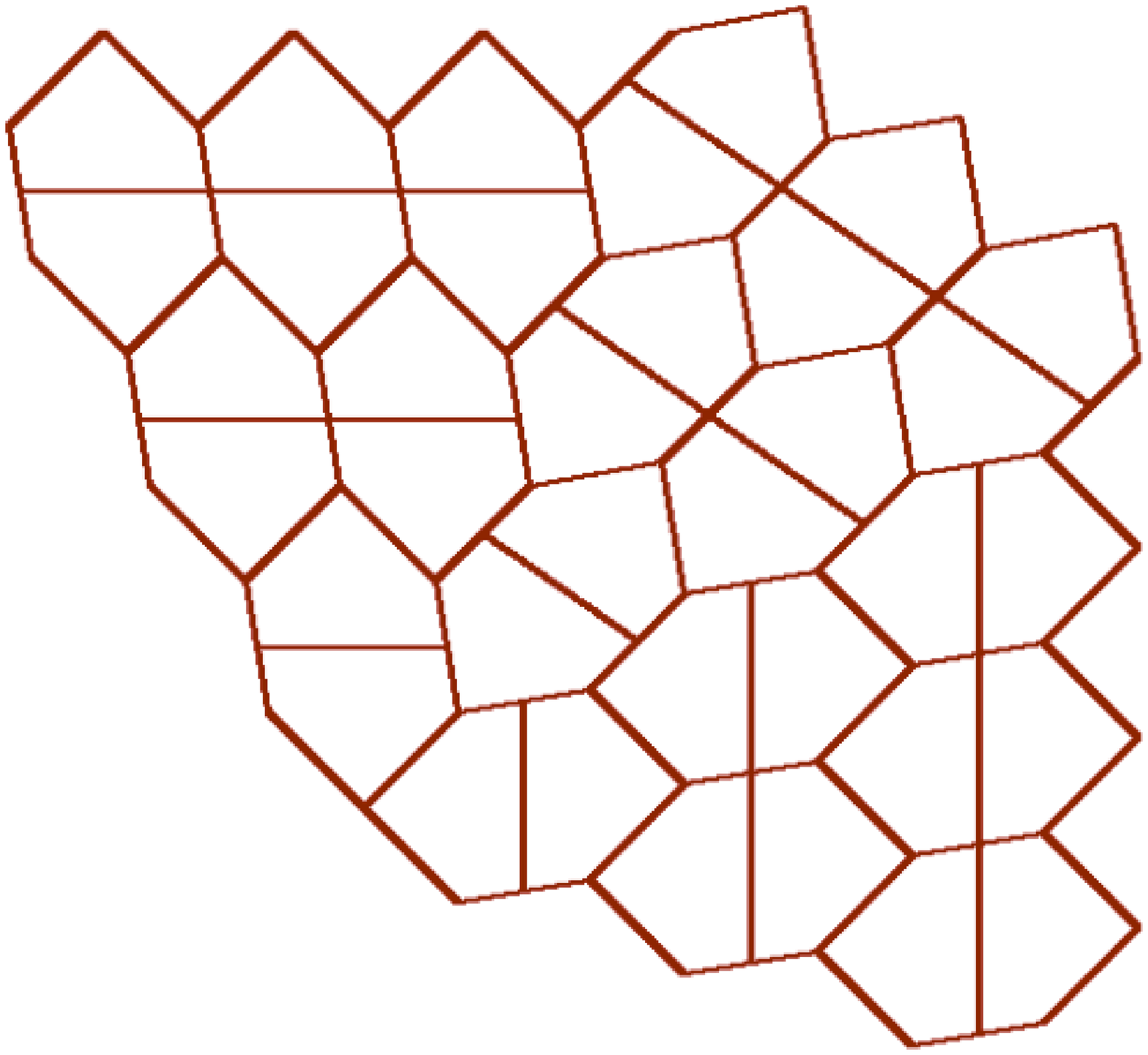} 
\begin{description}
\item [{Figure}] 2.3 \,\,\,Another sector copy fills the gap after reflection
plus rotation 
\end{description}
With these operations we have to generate $2n$ sectors to fill the
360\textdegree{} angle around the origin of the plane and the $n$
gaps. By construction, the symmetry $\mathbf{C}_{n}$ is obvious.
To introduce the additional reflection axes needed for $\mathbf{D}_{n}$\,,
we just have to view the case $C=180{^{\circ}}-B/2$ and $D=90{^{\circ}}$
(see Fig. 2.4). Here the sector itself has reflection symmetry. Finally,
we have $2n$ sector copies but each two of them share the same axis.
So there are $n$ axes for $\mathbf{D}_{n}$. (For a full view, see
also the examples in section 3.)

\,

\includegraphics[width=4cm]{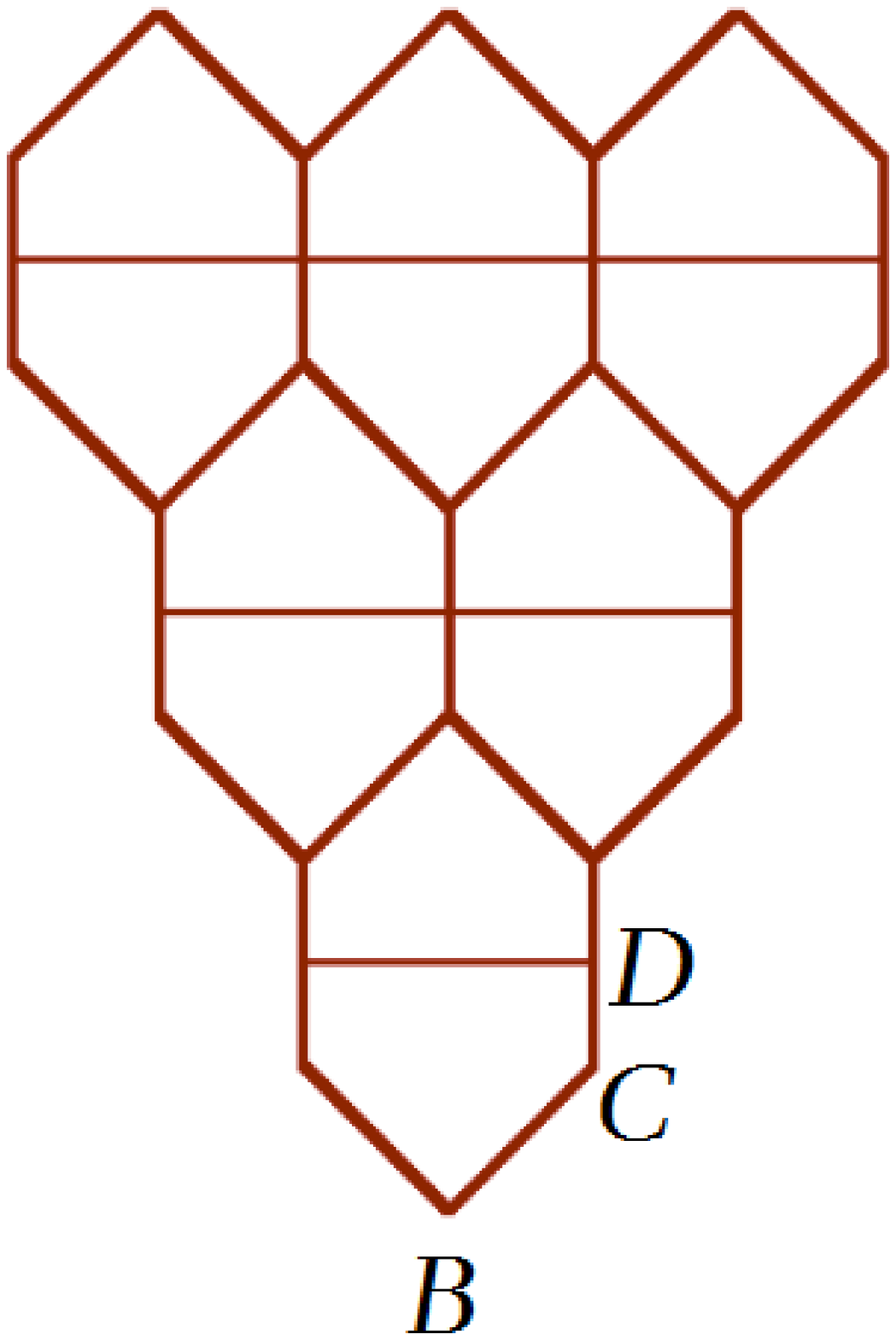} 
\begin{description}
\item [{Figure}] 2.4 \,\,\,Sector with reflection symmetry 
\end{description}
How can we see that our tilings are really plane filling? To prove
this, we can regard the tiling of each sector as a sequence of rows
with a growing number of hexagons. Consider those rows with a fixed
number of $m$ ($>1$) hexagons, which $-$ as a union $-$ always
forms a closed ring built out of $2nm$ hexagons. Then each of the
described tilings can be seen as a growing series of rings around
the centre with each ring consisting of hexagons. The first union
with $m=1$ is topologically equivalent to a disk. Each further ring
enlarges this disk, the size of which must grow with linearly increasing
diameter, since the shape of the hexagons cannot degenerate. Each
new ring fits to the previous one without gap. So, for any point on
the plane it will be covered by this growing disk after enough rings
were added consecutively.

But is it always possible to construct a pentagon with inner angles
$<180{^{\circ}}$? It is, since $120{^{\circ}}\ge B>0$ for all $n>2$.
Then $A$ and $C$ can be chosen e.\,g. $(360{^{\circ}}-B)/2\thinspace\pm\thinspace B/4$
respectively, which is smaller than $180{^\circ}$. The choice of
$D$ can be $270{^{\circ}}-B/2-C$, which is also between 0 and $180{^\circ}$.

Finally we have to discuss the cases $n=2$ and $n=1$. Here the above
construction will result in rectangles instead of hexagons, which
will not deliver proper pentagons with angles $<180{^{\circ}}$. So
we should look into the existing catalogue for pentagon tilings. The
most simple ones, called ``houses tiling``, will help. Simple examples
for $n=2$ are shown. The point of symmetry is marked by a dot.

\,

\includegraphics[width=9.5cm]{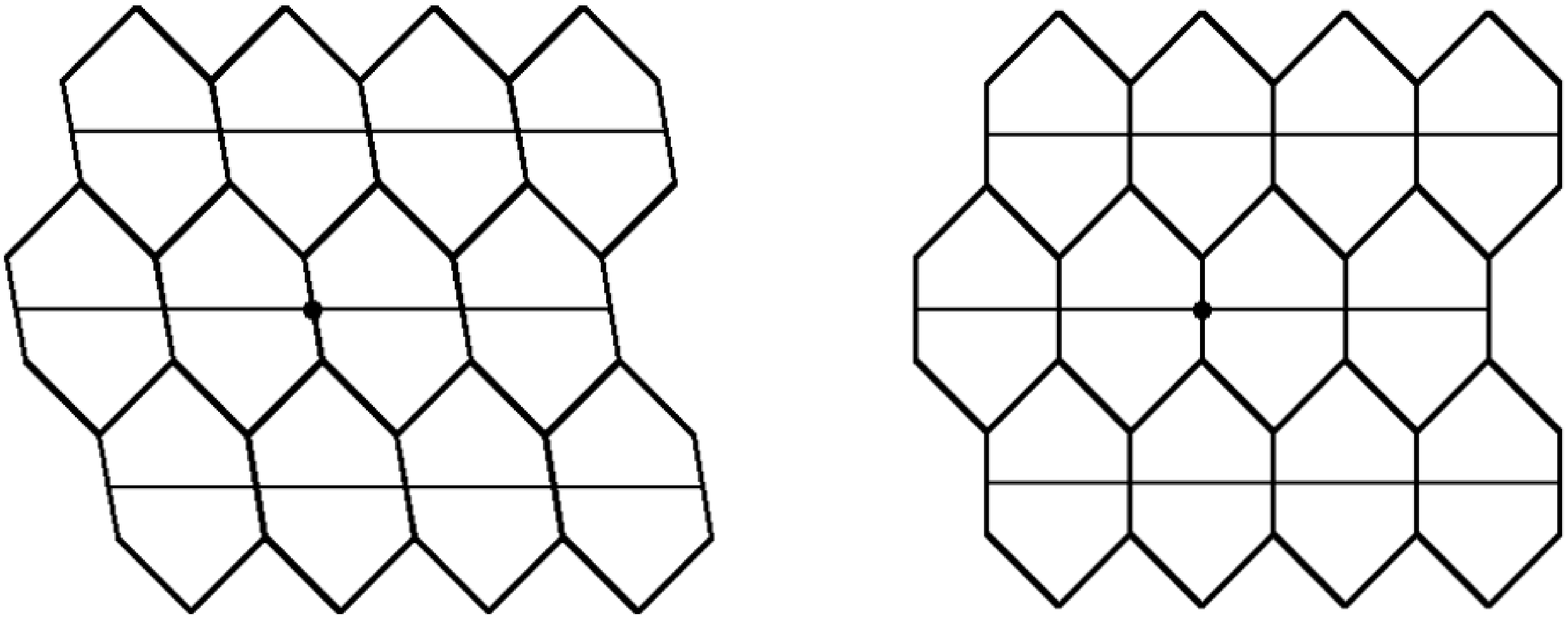} 
\begin{description}
\item [{Figure}] 2.5 \,\,\,Examples for $\mathbf{C}_{2}$ \,\,\,\,\,and
$\mathbf{D}_{2}$ symmetry 
\end{description}
The case $\mathbf{C}_{1}$ means ``neither rotational nor reflexion
symmetry\textquotedblright{} (see fig. 2.6 left) and $\mathbf{D}_{1}$
has only one symmetry axis and no rotational symmetry (see the right
part of fig. 2.6).

\,

\,

\includegraphics[width=10cm]{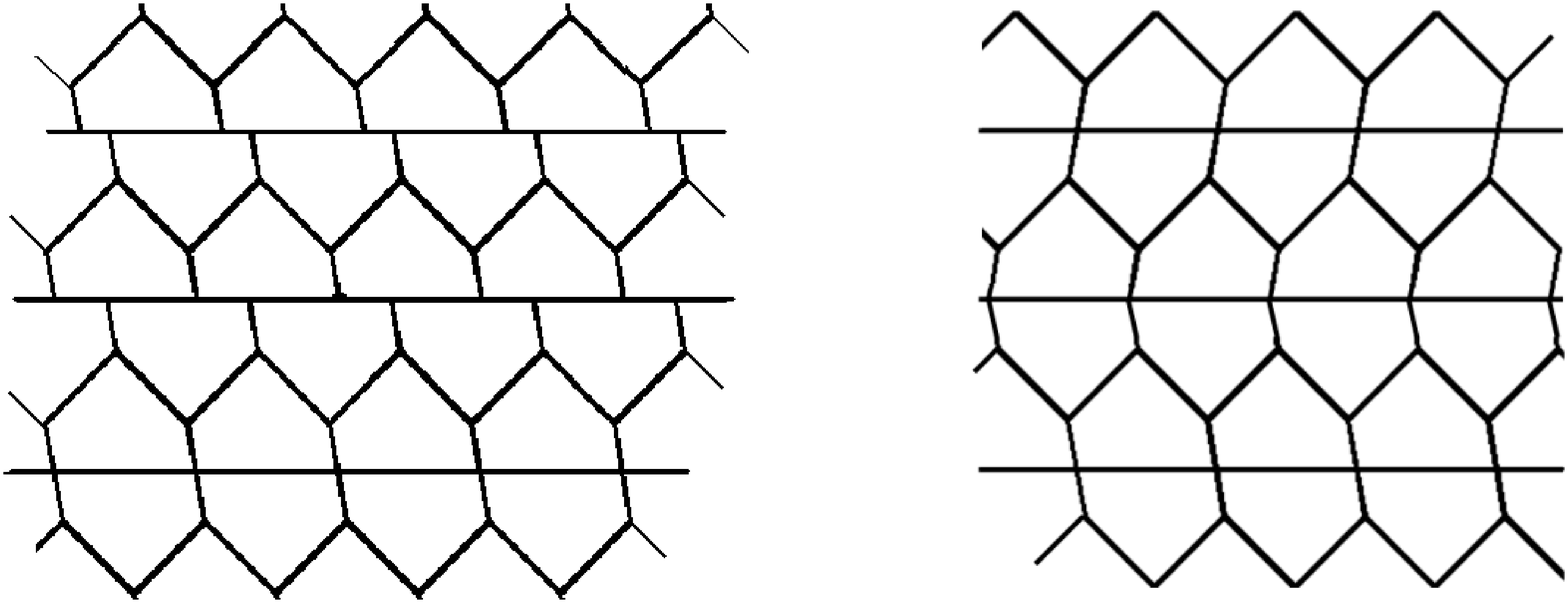} 
\begin{description}
\item [{Figure}] 2.6 \,\,\,Examples for $\mathbf{C}_{1}$ \,\,\,\,\,and
$\mathbf{D}_{1}$ symmetry 
\end{description}
So all symmetry types $\mathbf{C}_{n}$ and $\mathbf{D}_{n}$ were
constructed with pentagons having property 1, q.e.d.

\section{Examples}

We should not finish the paper without showing some of the nice tilings
resulting from the above construction. The first figure represents
$\mathbf{C}_{5}$ symmetry:

\,

\includegraphics[width=10cm]{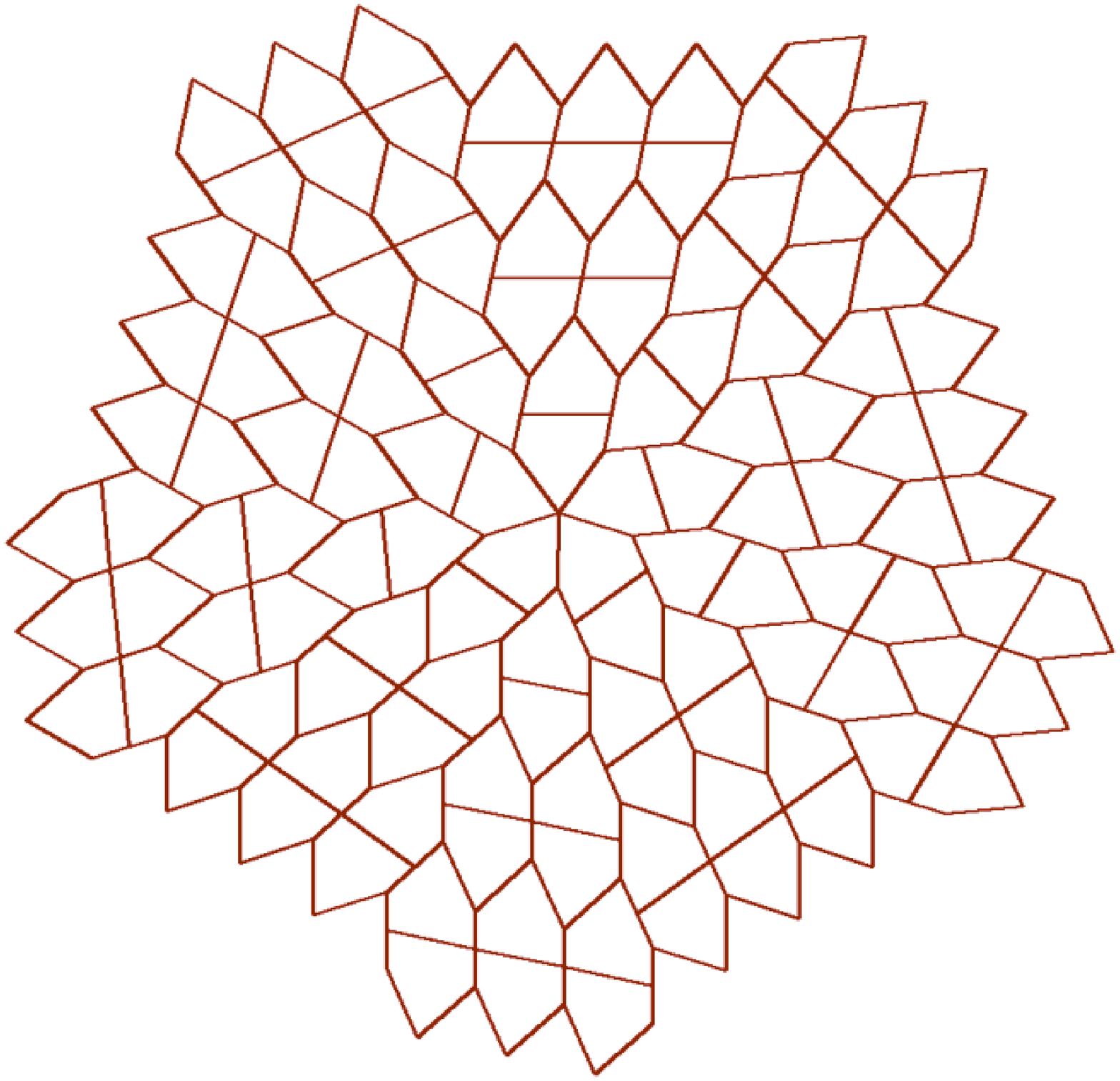} 
\begin{description}
\item [{Figure}] 3.1 \,\,\,$\mathbf{C}_{5}$ symmetry with $(A,\thinspace B,\thinspace C,\thinspace D,\thinspace E)=(132{^{\circ}},72{^{\circ}},156{^{\circ}},78{^{\circ}},102{^{\circ}})$ 
\end{description}
The second one is representing $\mathbf{D}_{7}$:

\,

\includegraphics[width=10.5cm]{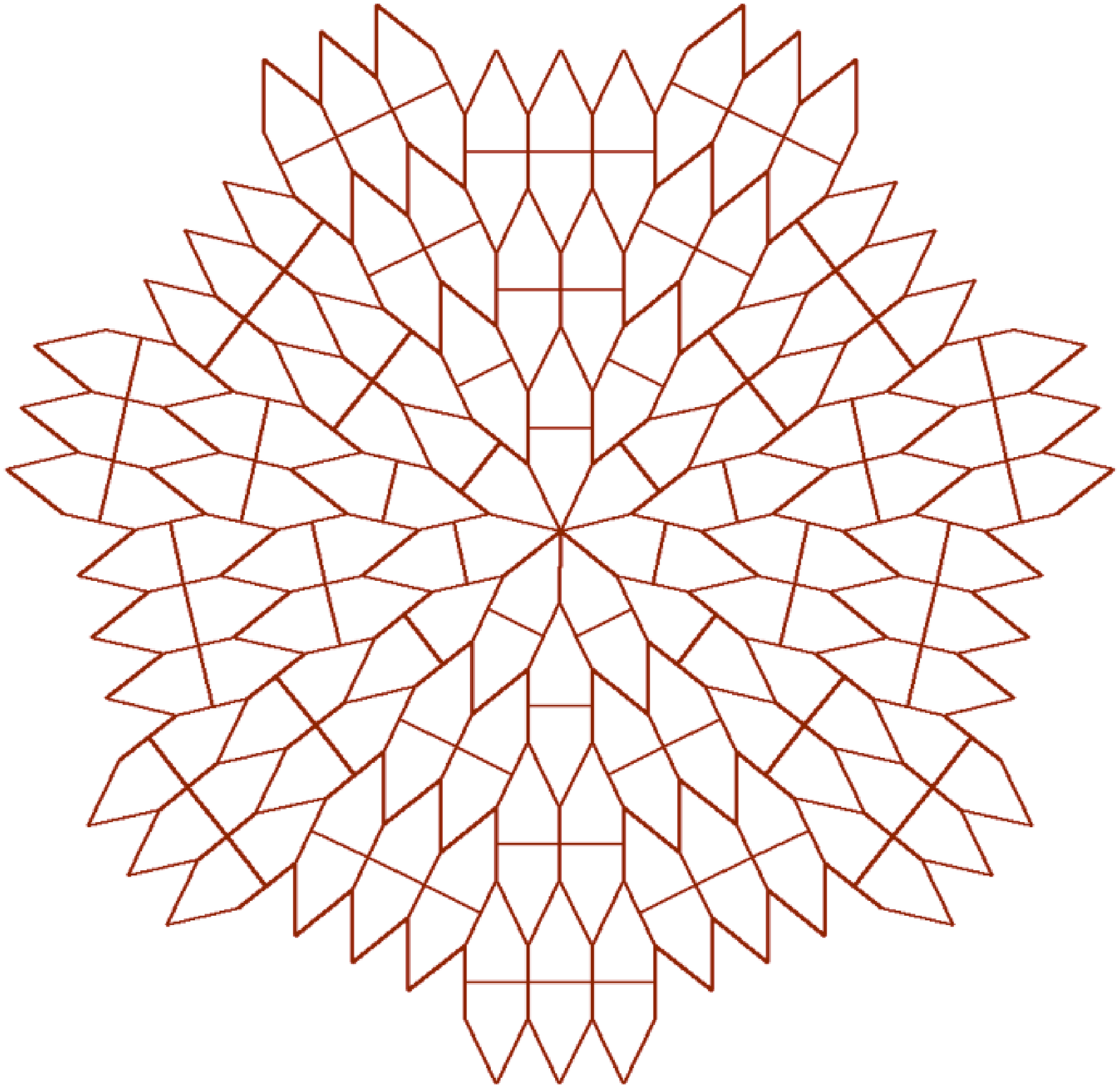} 
\begin{description}
\item [{Figure}] 3.2 \,\,\,Example for $\mathbf{D}_{7}$ symmetry 
\end{description}

\section{Spirals}

Another visually attractive property of the above construction is
the following. Any tiling according to chapter 2 with $n$-fold rotational
symmetry and $n>1$ can also be regarded as spiral tiling with $n$
congruent arms (see \cite{key-10} for a precise definition of the
term ``spiral tiling\textquotedblright ). To give an example, we
can identify one of the seven spiral arms within the $\mathbf{D}_{7}$
tiling from the above figure 3.2.

Take one of the innermost hexagons for the spiral's begin and walk
outward obeying the rule: Find the dividing line in the hexagon's
middle and choose the neighboring hexagon at the line's right endpoint
as next hexagon for the spiral (where ``right\textquotedblright{}
means ``right viewed from the origin\textquotedblright ).

Figure 4.1 displays this partition into spiral arms. In the well-known
book of Gr\"unbaum and Shephard \cite{key-7} it was put as an open
question if spiral tilings exist with $n$ arms for any odd $n>5$.
Later several spiral tilings have been published (e.g., \cite{key-8},
\cite{key-9} or \cite{key-3}) but \textendash{} to the author's
knowledge \textendash{} those with higher number of arms had non-convex
tiles. Here we can show spirals with arbitrary number of congruent
arms and with convex pentagons as prototiles.

\,

\includegraphics[width=10.5cm]{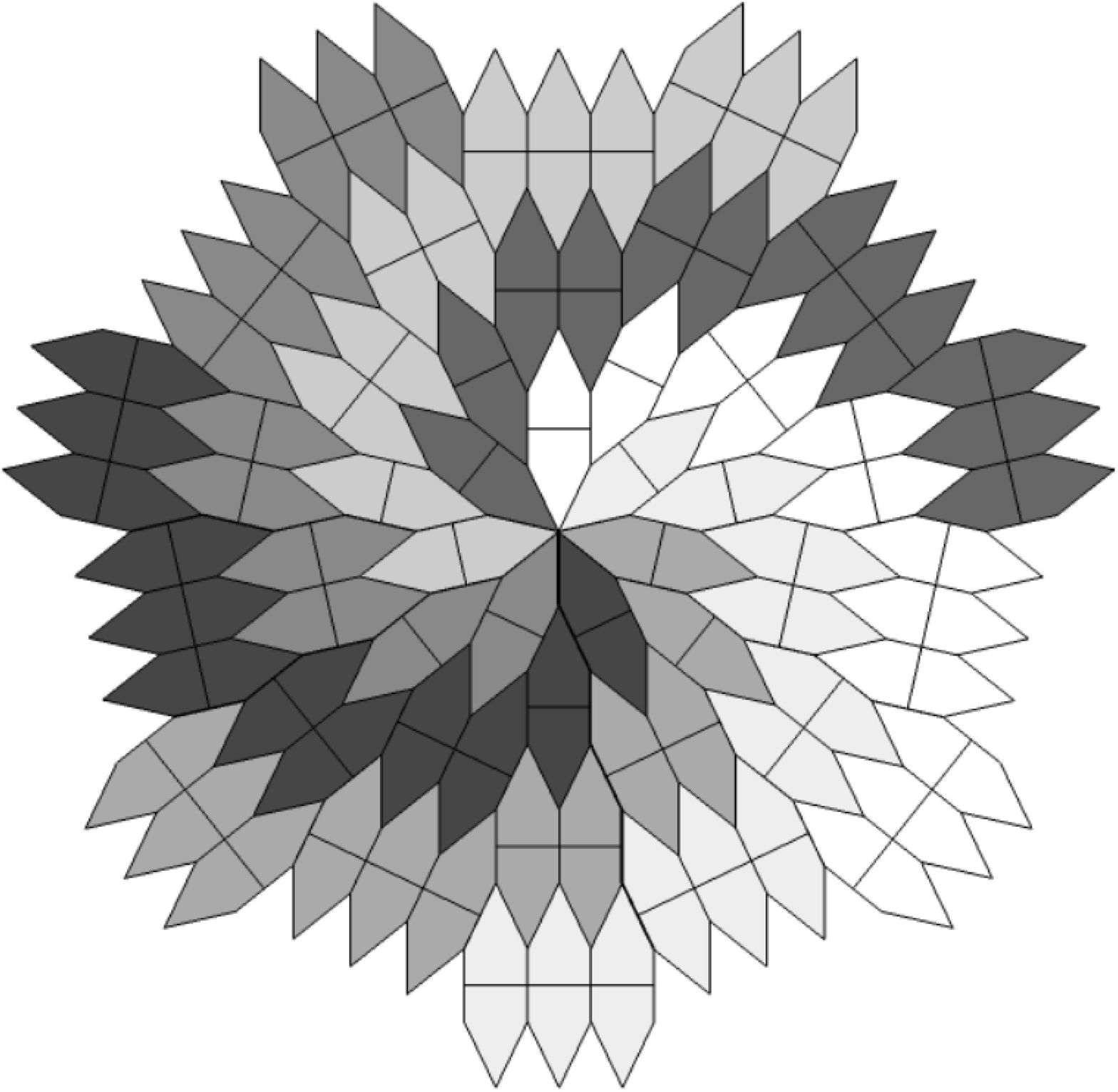} 
\begin{description}
\item [{Figure}] 4.1 \,\,\,Demonstration of the tiling's spiral character 
\end{description}
To summarize the results of this paper, we can state that for any
rotational symmetry type $\mathbf{C}_{n}$ or $\mathbf{D}_{n}$ there
is a tiling of the 2D plane with convex pentagons representing this
symmetry. Each of these tilings with $n>1$ can also be regarded as
spiral tiling with $n$ congruent arms. By construction it is obvious
that this holds for hexagons, as well.

After handling all rotational symmetric cases, there remains an open
question in the periodic case. The author's presumption is that it
will not be possible to represent all 2D symmetry types by convex
pentagon tilings, unless degenerated pentagons are used (i.e. triangles
or quadrangles).


\begin{thebibliography}{10}
\bibitem{key-1}https://en.wikipedia.org/wiki/Pentagonal\_tiling,
viewed Sept. 15th 2015

\bibitem{key-2}D. Schattschneider, ``Tiling the Plane with Congruent
Pentagons\textquotedblright , Mathematics Magazine 51(1) (1978), pp.
29-44

\bibitem{key-3}D.R. Simonds, ``Central Tessellations with an Equilateral
Pentagon\textquotedblright , Mathematics Teaching 81 (1977), pp. 36-37

\bibitem{key-4}M.D. Hirschhorn, D.C. Hunt, ``Equilateral Convex
Pentagons Which Tile the Plane\textquotedblright , Journal of Combinatorial
Theory, Series A, 39(1), (1985), pp. 1-18

\bibitem{key-5}https://en.wikipedia.org/wiki/Symmetry\#In\_geometry,
viewed Sept. 15th 2015

\bibitem{key-6}A. Bellos, ``Attack on the pentagon results in discovery
of new mathematical tile\textquotedblright , The Guardian, August
11th 2015

\bibitem{key-7}B. Gr\"unbaum, G.C. Shephard, ``Tilings and Patterns\textquotedblright ,
W.H. Freeman and Co., New York, 1987

\bibitem{key-8}P. Gailiunas, ``Spiral Tilings\textquotedblright ,
http://www.mi.sanu.ac.rs/vismath/gal/, viewed Sept. 15th 2015

\bibitem{key-9}D.L. Stock, B.A.Wichmann, ``Odd Spiral Tilings\textquotedblright ,
Mathematics Magazine 73 (5), Dec 2000, pp. 339-347

\bibitem{key-10}B. Klaassen, ``How to Define a Spiral Tiling?\textquotedblright ,
to appear 2016

\bibitem{key-11}G. Maloney, ``On substitution tilings of the plane
with n-fold rotational symmetry\textquotedblright , Discrete Math.
Theor. Computs. Sci. 17 (1), 2015, pp. 395-411

\bibitem{key-12}J. Garcia Escudero, ``Randomness and Topological
Invariants in Pentagonal Tiling Spaces\textquotedblright , Discrete
Dyn. Nat. Soc., Art. ID 946913, 23, 2011\end{thebibliography}
\end{document}